\documentclass{article}
\usepackage{latexsym}
\usepackage{amssymb}
\usepackage{amsmath}
\usepackage{amsfonts}
\newtheorem{teo}{Theorem}[section]
\newtheorem{lema}[teo]{Lemma}

\newtheorem{prop}[teo]{Proposition}
\newtheorem{rema}[teo]{Remark}

\newcommand{\cF}{{\cal{F}}}

\begin{document}
\title{A $d$-dimensional Brownian motion as a weak limit  from a one-dimensional Poisson process}
\date{}

\bigskip

\author{Xavier Bardina$^\dag$ , Carles Rovira \footnote{Corresponding author.} \smallskip
\\{\footnotesize $^\dag$ Dept. Matem\`atiques, Edifici C, Universitat Aut\`onoma de Barcelona, 08193-Bellaterra}
\\{\footnotesize $^*$ Facultat de Matem\`atiques, Universitat de Barcelona, Gran Via 585, 08007-Barcelona}
\\{\footnotesize E-mail: Xavier.Bardina@uab.cat, Carles.Rovira@ub.edu}}
\maketitle
\begin{abstract}
We show how from an unique standard Poisson process we can build a family of processes that converges in law to 
a $d$-dimensional standard Brownian motion for any $d \ge 1$.
\end{abstract}

\smallskip


\section{Introduction and main result}

Consider the sequences of processes
\begin{eqnarray*}
\{z_{\varepsilon}^{\theta}(t)=\varepsilon\int_{0}^{\frac{2t}{\varepsilon^2}}\cos(\theta
N_s)ds,\quad t\in[0,T]\},\\
\{y_{\varepsilon}^{\theta}(t)=\varepsilon\int_{0}^{\frac{2t}{\varepsilon^2}}\sin(\theta
N_s)ds,\quad t\in[0,T]\},
\end{eqnarray*}
where $\{N_s,\, s\geq 0\}$ is a standard Poisson process.

When $\theta=0$, the processes
$z_{\varepsilon}^{\theta}$  are deterministic and go to
infinity when  $\varepsilon$ tends to zero.
On the other hand, when $\theta=\pi$, the processes
$z_{\varepsilon}^{\theta}$ can be written as
$$z_{\varepsilon}^{\theta}(t)=\varepsilon\int_0^{\frac{2t}{\varepsilon^{2}}}(-1)^{N_s}ds.$$
This case was studied by Stroock in \cite{S}, where he proved that the laws
of these processes in the space of continuos functions on $[0,T]$
converge weakly towards the law of $\sqrt{2}W_{t}$, where $\{W_{t};\,
t\in[0,T]\}$ is a standard Brownian motion.
When $\theta=0$ or $\theta=\pi$, the processes $y_{\varepsilon}^{\theta}$ are constant and equal to zero.

When $\theta\in (0,\pi)\cup(\pi,2\pi)$ it is proved in \cite{Ba} that when $\varepsilon$ tends to zero
the processes $x_\varepsilon^\theta:=(z_\varepsilon^\theta,y_\varepsilon^\theta)$ converge weakly towards two independent standard Brownian motions.

The aim of this paper is to extend this result to a $d$-dimensional case for any $d\ge 1$. 

Now, given $\theta_1,\ldots,\theta_n,\theta_{n+1},\ldots,\theta_{n+m}$  let us consider the process:
$$
\{x_{\varepsilon}^{\theta_1,\ldots,\theta_n;\theta_{n+1},\ldots,\theta_{n+m}}(t)= 
(z_\varepsilon^{\theta_1},\ldots, z_\varepsilon^{\theta_n}, y_\varepsilon^{\theta_{n+1}},\ldots, y_\varepsilon^{\theta_{n+m}})(t),\quad t\in[0,T]\}.
$$

For simplicity we will denote by $\theta$ the $n+m$ values  $\theta_1,\ldots,\theta_n,\theta_{n+1},\ldots,\theta_{n+m}$. We will assume
the following hypothesis {\bf (H)} on $\theta$

\medskip
\begin{itemize}

\item  $\theta_i \in (0,\pi) \cup (\pi,2\pi), 1 \le i \le n+m$, 

\item   $\theta_i+\theta_j\neq 2\pi$ for all $1 \le i ,j \le n+m$. 

\item  $\theta_i-\theta_j\neq 0$ for all $1 \le i ,j \le n$ and $n+1 \le i ,j \le n+m$. 

\end{itemize} 

Let us point out the meaning of the last hypothesis: two parameters $\theta_i$ and $\theta_j$ can only be equal if $i \le n$ and $j \ge n+1$. In other words we can deal with $\varepsilon\int_{0}^{\frac{2t}{\varepsilon^2}}\cos(\theta_iN_s)ds$ and $\varepsilon\int_{0}^{\frac{2t}{\varepsilon^2}}\sin(\theta_i
N_s)ds$, but the results will not be possible (obviously) if we have two times $\varepsilon\int_{0}^{\frac{2t}{\varepsilon^2}}\cos(\theta_iN_s)ds$ or $\varepsilon\int_{0}^{\frac{2t}{\varepsilon^2}}\sin(\theta_i
N_s)ds$.

\medskip

 Our result states as follows,
\begin{teo}\label{teo1}

Consider $P_{\varepsilon}^{\theta}$ the image law of
$x_{\varepsilon}^{\theta}$ in the Banach space \linebreak $\mathcal C([0,T],\mathbb
R^{n+m})$ of continuous functions on
$[0,T]$. If $\theta$ satisfies hypothesis (H) then
$P_\varepsilon^{\theta}$ converges weakly as $\varepsilon$ tends to
zero towards the law on
$\mathcal C([0,T],\mathbb R^{n+m})$ of a $n+m$-dimensional standard Brownian motion.

\end{teo}

\medskip

\begin{rema}
It will be also possible to consider the case $\theta_i = \pi$ for some $i \in \{1,\ldots,n\}$.
In this case, we need to deal with $\frac{1}{\sqrt 2} z_\varepsilon^{\theta_i}$ instead of 
$z_\varepsilon^{\theta_i}$. Nevertheless, the proof follows the same computations.
\end{rema}

\smallskip

The structure of the paper is the following. In Section 2 we recall the basic results of \cite{Ba} that implies
our theorem when $m=n=1$ and $\theta_1=\theta_2$. In Section 3 we give the proof of our main theorem.

Along the paper $K$ denote positive constants, not depending on
$\varepsilon$,
which may change from one expression to another one.

\section{The two-dimensional case}

In \cite{Ba} it is proved  an approximation in law of the complex Brownian
motion by processes constructed from a unique standard Poisson process.

\begin{teo} \cite[Theorem 1.1]{Ba}
Define for any $\varepsilon >0$
$$
\{v_{\varepsilon}^{\theta}(t)=\varepsilon\int_{0}^{\frac{2t}{\varepsilon^2}}e^{i\theta
N_s}ds,\quad t\in[0,T]\}$$
where $\{N_s,\, s\geq 0\}$ is a standard Poisson process.
Consider $P_{\varepsilon}^{\theta}$ the image law of
$v_{\varepsilon}^{\theta}$ in the Banach space $\mathcal C([0,T],\mathbb
C)$ of continuous functions on
$[0,T]$. Then, if $\theta \in (0,\pi) \cup (\pi, 2 \pi)$, 
$P_\varepsilon^{\theta}$ converges weakly as $\varepsilon$ tends to
zero towards the law on
$\mathcal C([0,T],\mathbb C)$ of a complex Brownian motion.

\end{teo}

In fact, it corresponds to our two-dimensional case with $n=m=1$ and $\theta_1=\theta_2=\theta$ for $\theta\in(0,\pi)\cup(\pi,2\pi)$.
Set $P_{\varepsilon}^{\theta}$ the image law of
$x_{\varepsilon}^{\theta}:=(z_{\varepsilon}^{\theta},y_{\varepsilon}^{\theta})$ in the space $\mathcal C([0,T],\mathbb
R^{2}).$

The proof of the weak convergence is obtained checking that the family
$P_{\varepsilon}^{\theta}$ is tight and that the law of
all possible weak limits of $P_\varepsilon^{\theta}$ is the law of two independent
standard Brownian motions.

The tigthness is proved using the Billingsley criterium (see Theorem 12.3 of
\cite{B}). Since  the processes are null on the origin it suffices
to prove the following lemma (see  Lemma 2.1 in \cite{Ba}).
\begin{lema} \label{lema2.1} There exists a constant $K$ such that for any
$s<t$
$$\sup_{\varepsilon}\big(E(\varepsilon\int_{\frac{2s}{\varepsilon^2}}^{\frac{2t}{\varepsilon^2}}\cos(\theta
N_x)dx)^{4} +
E(\varepsilon\int_{\frac{2s}{\varepsilon^2}}^{\frac{2t}{\varepsilon^2}}\sin(\theta
N_x)dx)^{4}\big)\leq K(t-s)^{2}.$$
\end{lema}

In order to identify the limit law, it is considered $\{P_{\varepsilon_n}^{\theta}\}_{n}$ a subsequence of
$\{P_{\varepsilon}^{\theta}\}_{\varepsilon}$ (that is also denoted by
$\{P_{\varepsilon}^{\theta}\}$) weakly convergent to some probability
$P^{\theta}$. Then, it is checked that  the two components of the
canonical process $X=(Z,Y)=\{X_{t}(x)=x(t)=(z(t),y(t))\}$ 
under the probability $P^{\theta}$ are two independent Brownian motions.

Using Paul L\'evy's theorem (see Theorem \ref{levy} below) it suffices to prove that under $P^{\theta}$,
$Z$ and $Y$ are both martingales with respect
to the natural filtration, $\{\mathcal F_{t}\}$, with quadratic
variations $<Z,Z>_{t}=t$, $<Y,Y>_{t}=t$ and
covariation $<Z,Y>_{t}=0$.

To check the martingale property with respect to the natural filtration $\{\mathcal
F_{t}\}$, it is proved (see  subsection 3.1 in \cite{Ba}) that for any $s_{1}\leq s_{2}\leq\cdots\leq s_{k}
\leq s<t$ and for any bounded continuous function $\varphi:\mathbb R^{2k}
\longrightarrow\mathbb R$,
$$E_{P^{\theta}}\big[\varphi(X_{s_{1}},...,X_{s_{k}})(Z_{t}-Z_{s})\big]=0,$$
$$E_{P^{\theta}}\big[\varphi(X_{s_{1}},...,X_{s_{k}})(Y_{t}-Y_{s})\big]=0.$$

The computation of the quadratic variations and covariation is done in the following proposition
(see  Proposition 3.1 in \cite{Ba}).

\begin{prop} \label{prop3.1}
Consider $\{P_{\varepsilon}^{\theta}\}$ the laws on $\mathcal
C([0,T],\mathbb R^2)$ of the processes $x_{\varepsilon}^{\theta}$ and assume that $P_{\varepsilon_{n}}^{\theta}$ is a
subsequence weakly convergent to $P^{\theta}$. Let $X=(Z,Y)$ be the canonical
process and let
$\{\mathcal F_{t}\}$ be its natural filtration. Then, under $P^{\theta}$,
if
$\theta\in(0,\pi)\cup(\pi,2\pi)$  it is hold that the quadratic variations
$<Z,Z>_{t}=t$, $<Y,Y>_{t}=t$ and the covariation
$<Z,Y>_{t}=0$.
\end{prop}

\section{Proof of the main result}

In this section we will give the proof of Theorem \ref{teo1}. We will follow the same method than in 
\cite{Ba}. So, it suffices to check the tightness
of the family
$P_{\varepsilon}^{\theta}$  and to identify  the law of
all possible weak limits of $P_\varepsilon^{\theta}$.

The tigthness is proved also using the Billingsley criterium, and it is an obvious consequence of
Lemma \ref{lema2.1}. 

The identification of the limit law will be done using Paul L\'evy's theorem.

\begin{teo} \cite{L} \label{levy}
Let $X=\{X_t=(X_t^{(1)},\ldots,X_t^{(d)}\},\cF_t, 0 \le t < \infty \}$ be a continuous, adapted process
in $\mathbb R^d$ such that, for any component $1 \le k \le d$ the process
$$M_t^{(k)}:= X_t^{(k)}-X_0^{(k)}, \qquad 0 \le t < \infty$$
is a continuous local martingale relative to $\{ \cF_t  \}$ and the cross-variations are given by
$<M^{(k)},M^{(j)}>_{t}= \delta_{k,j} t,$ for $1 \le k,j \le d$.
Then $\{ X_t, \cF_t, 0 \le t < \infty \}$ is a $d$-dimensional Brownian motion.

\end{teo}

Let us consider $\{P_{\varepsilon_n}^{\theta}\}_{n}$ a subsequence of
$\{P_{\varepsilon}^{\theta}\}_{\varepsilon}$ (that we also will denote by
$\{P_{\varepsilon}^{\theta}\}$) weakly convergent to some probability $P^{\theta}$. Consider the
canonical process $X=(Z^1,\ldots,Z^n,Y^{n+1},\ldots,Y^{n+m})$ 
under the probability $P^{\theta}$. It suffices to check that under $P^{\theta}$,
$Z^i, 1 \le i \le n,$ and $Y^j, n+1 \le j \le n+m,$ are  martingales with respect
to the natural filtration, $\{\mathcal F_{t}\}$, with quadratic
variations $<Z^i,Z^i>_{t}=t, 1 \le i \le n$, $<Y^j,Y^j>_{t}=t, n+1 \le j \le n+m,$ and
covariation $<Z^i,Z^l>_{t}=<Y^j,Y^h>_{t}=<Z^i,Y^j>_{t}=0, 1 \le i \not= l \le n, n+1 \le j \not= h\le n+m$.

In order to check the martingale property with respect to the natural filtration $\{\mathcal
F_{t}\}$  it suffices to prove that for any $s_{1}\leq s_{2}\leq\cdots\leq s_{k}
\leq s < t$ and for any bounded continuous function $\varphi:\mathbb R^{(n+m)k}
\longrightarrow\mathbb R$,
$$E_{P^{\theta}}\big[\varphi(X_{s_{1}},...,X_{s_{k}})(Z^i_{t}-Z^i_{s})\big]=0, \quad{\mbox {\rm for} }\quad 1 \le i \le n,$$
$$E_{P^{\theta}}\big[\varphi(X_{s_{1}},...,X_{s_{k}})(Y^j_{t}-Y^j_{s})\big]=0, \quad{\mbox {\rm for} }\quad n+1 \le j \le n+m.$$
These computations has been done in  subsection 3.1 in \cite{Ba}.

The proof of the quadratic variations can be done following exactly the proof of Proposition 3.1 in \cite{Ba}.
So, it remains only to compute all the covariations. 

First we have to prove that $<Z^i,Z^l>_{t}=0$, for $1 \le i \not= l \le n$. It suffices to prove
that for any     $s_{1}\leq s_{2}\leq\cdots\leq s_{k}
\leq s<t$ and for any bounded continuous function $\varphi:\mathbb R^{(n+m)k}
\longrightarrow\mathbb R$,
$$E\big[\varphi(x_\varepsilon^{\theta}(s_{1})...,x_\varepsilon^{\theta}(s_{k}))
(z_\varepsilon^{\theta_i}(t)-z_\varepsilon^{\theta_i}(s))
(z_\varepsilon^{\theta_l}(t)-z_\varepsilon^{\theta_l}(s))\big],$$
converges to zero when $\varepsilon$ tends to zero. Notice that this last expression
can be written as 
\begin{equation}\label{eqcos}
E\left(\varphi\left(x_{\varepsilon}^{\theta}(s_1),\dots,x_{\varepsilon}^{\theta}(s_k)\right)
\left(\varepsilon\int_{\frac{2s}{\varepsilon^2}}^{\frac{2t}{\varepsilon^2}}\cos(\theta_iN_x)dx\right)
\left(\varepsilon\int_{\frac{2s}{\varepsilon^2}}^{\frac{2t}{\varepsilon^2}}\cos(\theta_lN_x)dx\right)\right).
\end{equation}

Similarly, to prove that $<Y^j,Y^h>_{t}=0$, for $n+1 \le j \not= h\le n+m$, it is enough to prove
that for any     $s_{1}\leq s_{2}\leq\cdots\leq s_{k}
\leq s < t$ and for any bounded continuous function $\varphi:\mathbb R^{(n+m)k}
\longrightarrow\mathbb R$,
\begin{equation}\label{eqsin}
E\left(\varphi\left(x_{\varepsilon}^{\theta}(s_1),\dots,x_{\varepsilon}^{\theta}(s_k)\right)
\left(\varepsilon\int_{\frac{2s}{\varepsilon^2}}^{\frac{2t}{\varepsilon^2}}\sin(\theta_jN_x)dx\right)
\left(\varepsilon\int_{\frac{2s}{\varepsilon^2}}^{\frac{2t}{\varepsilon^2}}\sin(\theta_hN_x)dx\right)\right)
\end{equation}
converges to zero when $\varepsilon$ tends to zero.

Finally, to prove that $<Z^i,Y^j>_{t}=0,$ for $ 1 \le i \le n <n+1 \le j \le n+m$, it is enough to prove
that for any     $s_{1}\leq s_{2}\leq\cdots\leq s_{k}
\leq s<t$ and for any bounded continuous function $\varphi:\mathbb R^{(n+m)k}
\longrightarrow\mathbb R$,
\begin{equation}\label{eqcossin}
E\left(\varphi\left(x_{\varepsilon}^{\theta}(s_1),\dots,x_{\varepsilon}^{\theta}(s_k)\right)
\left(\varepsilon\int_{\frac{2s}{\varepsilon^2}}^{\frac{2t}{\varepsilon^2}}\cos(\theta_iN_x)dx\right)
\left(\varepsilon\int_{\frac{2s}{\varepsilon^2}}^{\frac{2t}{\varepsilon^2}}\sin(\theta_jN_x)dx\right)\right)
\end{equation}
converges to zero when $\varepsilon$ tends to zero.

Let us finish the proof of the theorem checking the convergence to zero of (\ref{eqcos}), (\ref{eqsin}) and (\ref{eqcossin})
when $\varepsilon$ tends to zero. For simplicity we will use only $\theta_1$ and $\theta_2$.

\bigskip

{\it Study of (\ref{eqsin}).}
Notice that (\ref{eqsin}) is equal to
\begin{eqnarray*}
&&\varepsilon^2\int_{\frac{2s}{\varepsilon^2}}^{\frac{2t}{\varepsilon^2}}\int_{\frac{2s}{\varepsilon^2}}^{x_2}
E\left(\varphi\left(x_{\varepsilon}^{\theta}(s_1),\dots,x_{\varepsilon}^{\theta}(s_k)\right)
\sin(\theta_1N_{x_1})\sin(\theta_2N_{x_2})\right)dx_1dx_2\\
&&+\varepsilon^2\int_{\frac{2s}{\varepsilon^2}}^{\frac{2t}{\varepsilon^2}}\int_{\frac{2s}{\varepsilon^2}}^{x_1}
E\left(\varphi\left(x_{\varepsilon}^{\theta}(s_1),\dots,x_{\varepsilon}^{\theta}(s_k)\right)
\sin(\theta_1N_{x_1})\sin(\theta_2N_{x_2})\right)dx_2dx_1\\
&:=&I_1+I_2.
\end{eqnarray*}

Using that $\sin(a)\sin(b)=\frac{\cos(a-b)-\cos(a+b)}{2}$ we obtain
that
\begin{eqnarray*}
I_1&=&\frac12\varepsilon^2\int_{\frac{2s}{\varepsilon^2}}^{\frac{2t}{\varepsilon^2}}\int_{\frac{2s}{\varepsilon^2}}^{x_2}
E\left(\varphi\left(x_{\varepsilon}^{\theta}(s_1),\dots,x_{\varepsilon}^{\theta}(s_k)\right)
\cos(\theta_1N_{x_1}-\theta_2N_{x_2})\right)dx_1dx_2\\
&&-\frac12\varepsilon^2\int_{\frac{2s}{\varepsilon^2}}^{\frac{2t}{\varepsilon^2}}\int_{\frac{2s}{\varepsilon^2}}^{x_2}
E\left(\varphi\left(x_{\varepsilon}^{\theta}(s_1),\dots,x_{\varepsilon}^{\theta}(s_k)\right)
\cos(\theta_1N_{x_1}+\theta_2N_{x_2})\right)dx_1dx_2\\
&:=&I_{1,1}-I_{1,2}.
\end{eqnarray*}

We start with the term $I_{1,1}$. Notice that
\begin{eqnarray*}
&&I_{1,1}\\&=&\frac12\textrm{Re}\left[\varepsilon^2\int_{\frac{2s}{\varepsilon^2}}^{\frac{2t}{\varepsilon^2}}\int_{\frac{2s}{\varepsilon^2}}^{x_2}
E\left(\varphi\left(x_{\varepsilon}^{\theta}(s_1),\dots,x_{\varepsilon}^{\theta}(s_k)\right)
e^{i(\theta_1N_{x_1}-\theta_2N_{x_2})}\right)dx_1dx_2\right]\\
&=&\frac12\textrm{Re}\left[\varepsilon^2\int_{\frac{2s}{\varepsilon^2}}^{\frac{2t}{\varepsilon^2}}\int_{\frac{2s}{\varepsilon^2}}^{x_2}
E\left(\varphi\Big(x_{\varepsilon}^{\theta}(s_1),\dots,x_{\varepsilon}^{\theta}(s_k)\right)
\right.\\&&\qquad\times\left.\left.e^{i(\theta_1-\theta_2)N_{\frac{2s}{\varepsilon^2}}}e^{i(\theta_1-\theta_2)(N_{x_1}-N_{\frac{2s}{\varepsilon^2}})}
e^{-i\theta_2(N_{x_2}-N_{x_1})}\right)dx_1dx_2\right].
\end{eqnarray*}
Using  the independence of the increments of the Poisson process and that
$E(e^{i \theta N_s})=e^{- s(1-e^{i \theta})}$ we get 
\begin{eqnarray*}
I_{1,1} &=&\frac12\textrm{Re}\left[\varepsilon^2\int_{\frac{2s}{\varepsilon^2}}^{\frac{2t}{\varepsilon^2}}\int_{\frac{2s}{\varepsilon^2}}^{x_2}
E\left(\varphi\left(x_{\varepsilon}^{\theta}(s_1),\dots,x_{\varepsilon}^{\theta}(s_k)\right)
e^{i(\theta_1-\theta_2)N_{\frac{2s}{\varepsilon^2}}}
\right)\right.\\&&\qquad\times\left.E\left(e^{i(\theta_1-\theta_2)(N_{x_1}-N_{\frac{2s}{\varepsilon^2}})}\right)
E\left(e^{-i\theta_2(N_{x_2}-N_{x_1})}\right)dx_1dx_2\right]\\
&\leq&K\varepsilon^2\int_{\frac{2s}{\varepsilon^2}}^{\frac{2t}{\varepsilon^2}}\int_{\frac{2s}{\varepsilon^2}}^{x_2}
\|e^{-\left(x_1-\frac{2s}{\varepsilon^2}\right)(1-e^{i(\theta_1-\theta_2)})}\|
\|e^{-\left(x_2-x_1\right)(1-e^{-i\theta_2})}\|
dx_1dx_2\\
&=&K\varepsilon^2\int_{\frac{2s}{\varepsilon^2}}^{\frac{2t}{\varepsilon^2}}\int_{\frac{2s}{\varepsilon^2}}^{x_2}
e^{-\left(x_1-\frac{2s}{\varepsilon^2}\right)(1-\cos(\theta_1-\theta_2))}
e^{-\left(x_2-x_1\right)(1-\cos(\theta_2))}
dx_1dx_2\\
&\leq&K\varepsilon^2\frac1{1-\cos(\theta_2)}\int_{\frac{2s}{\varepsilon^2}}^{\frac{2t}{\varepsilon^2}}
e^{-\left(x_1-\frac{2s}{\varepsilon^2}\right)(1-\cos(\theta_1-\theta_2))}dx_1\\
&\leq&K\varepsilon^2\frac1{1-\cos(\theta_2)}\frac1{1-\cos(\theta_1-\theta_2)},
\end{eqnarray*}
that clearly converges to zero when $\varepsilon$ tends to zero
because $\theta_2\neq 0$ and $\theta_1-\theta_2\neq 0$.

Using the decomposition
$$
e^{i ( \theta_1 N_{x_1} + \theta_2 N_{x_2} )} =  e^{i \theta_2 (N_{x_2}- N_{x_1} )} e^{i ( \theta_1+ \theta_2)( N_{x_1} - N_{\frac{2s}{\varepsilon^2}} )}e^{i ( \theta_1 + \theta_2) N_{\frac{2s}{\varepsilon^2}}}$$
and following the same computations we obtain that
$$I_{1,2}\leq K\varepsilon^2\frac1{1-\cos(\theta_2)}\frac1{1-\cos(\theta_1+\theta_2)}.$$
This last expression also goes to zero because $\theta_2\in (0,\pi) \cup (\pi, 2\pi)$
and $\theta_1+\theta_2\neq 2\pi$.

On the other hand, the expression $I_2$ is equal than the expression
$I_1$ interchanging the roles of $\theta_1$ and $\theta_2$. So, we
obtain that,
$$I_{2}\leq K\varepsilon^2\frac1{1-\cos(\theta_1)}\left(\frac1{1-\cos(\theta_2-\theta_1)}+
\frac1{1-\cos(\theta_1+\theta_2)}\right).$$ This last expression
also goes to zero because $\theta_1\in (0,\pi) \cup (\pi, 2 \pi)$,
$\theta_2-\theta_1\neq 0$ and $\theta_1+\theta_2\neq 2\pi$.

\bigskip

\pagebreak

{\it Study of (\ref{eqcos}).}
By the same computations and using that $\cos(a)\cos(b)=\frac{\cos(a+b)+\cos(a-b)}{2}$
we have also that,
\begin{eqnarray*}
&&E\left(\varphi\left(x_{\varepsilon}^{\theta}(s_1),\dots,x_{\varepsilon}^{\theta}(s_k)\right)
\left(\varepsilon\int_{\frac{2s}{\varepsilon^2}}^{\frac{2t}{\varepsilon^2}}\cos(\theta_1N_x)dx\right)
\left(\varepsilon\int_{\frac{2s}{\varepsilon^2}}^{\frac{2t}{\varepsilon^2}}\cos(\theta_2N_x)dx\right)\right)\\
&=&\varepsilon^2\int_{\frac{2s}{\varepsilon^2}}^{\frac{2t}{\varepsilon^2}}\int_{\frac{2s}{\varepsilon^2}}^{x_2}
E\left(\varphi\left(x_{\varepsilon}^{\theta}(s_1),\dots,x_{\varepsilon}^{\theta}(s_k)\right)
\cos(\theta_1N_{x_1})\cos(\theta_2N_{x_2})\right)dx_1dx_2\\
&&+\varepsilon^2\int_{\frac{2s}{\varepsilon^2}}^{\frac{2t}{\varepsilon^2}}\int_{\frac{2s}{\varepsilon^2}}^{x_1}
E\left(\varphi\left(x_{\varepsilon}^{\theta}(s_1),\dots,x_{\varepsilon}^{\theta}(s_k)\right)
\cos(\theta_1N_{x_1})\cos(\theta_2N_{x_2})\right)dx_2dx_1\\
&=&I_{1,1}+I_{1,2}+I_{2,1}+I_{2,2},
\end{eqnarray*}
where $I_{1,1}$ and $I_{1,2}$ has been defined in the study of (\ref{eqsin}) and $I_{2,1}$ and $I_{2,2}$ will be the anologous decomposition of $I_2$. So, we  have the convergence
to zero when $\varepsilon$ tends to zero.

\bigskip

{\it Study of (\ref{eqcossin}).}
Again the same computations yield that,
\begin{eqnarray*}
&&E\left(\varphi\left(x_{\varepsilon}^{\theta}(s_1),\dots,x_{\varepsilon}^{\theta}(s_k)\right)
\left(\varepsilon\int_{\frac{2s}{\varepsilon^2}}^{\frac{2t}{\varepsilon^2}}\sin(\theta_1N_x)dx\right)
\left(\varepsilon\int_{\frac{2s}{\varepsilon^2}}^{\frac{2t}{\varepsilon^2}}\cos(\theta_2N_x)dx\right)\right)\\
&=&\varepsilon^2\int_{\frac{2s}{\varepsilon^2}}^{\frac{2t}{\varepsilon^2}}\int_{\frac{2s}{\varepsilon^2}}^{x_2}
E\left(\varphi\left(x_{\varepsilon}^{\theta}(s_1),\dots,x_{\varepsilon}^{\theta}(s_k)\right)
\sin(\theta_1N_{x_1})\cos(\theta_2N_{x_2})\right)dx_1dx_2\\
&&+\varepsilon^2\int_{\frac{2s}{\varepsilon^2}}^{\frac{2t}{\varepsilon^2}}\int_{\frac{2s}{\varepsilon^2}}^{x_1}
E\left(\varphi\left(x_{\varepsilon}^{\theta}(s_1),\dots,x_{\varepsilon}^{\theta}(s_k)\right)
\sin(\theta_1N_{x_1})\cos(\theta_2N_{x_2})\right)dx_2dx_1\\
&=&J_1+J_2.
\end{eqnarray*}

Using that $\sin(a)\cos(b)=\frac{\sin(a+b)+\sin(a-b)}{2}$ we obtain
that
\begin{eqnarray*}
J_1&=&\frac12\varepsilon^2\int_{\frac{2s}{\varepsilon^2}}^{\frac{2t}{\varepsilon^2}}\int_{\frac{2s}{\varepsilon^2}}^{x_2}
E\left(\varphi\left(x_{\varepsilon}^{\theta}(s_1),\dots,x_{\varepsilon}^{\theta}(s_k)\right)
\sin(\theta_1N_{x_1}-\theta_2N_{x_2})\right)dx_1dx_2\\
&&+\frac12\varepsilon^2\int_{\frac{2s}{\varepsilon^2}}^{\frac{2t}{\varepsilon^2}}\int_{\frac{2s}{\varepsilon^2}}^{x_2}
E\left(\varphi\left(x_{\varepsilon}^{\theta}(s_1),\dots,x_{\varepsilon}^{\theta}(s_k)\right)
\sin(\theta_1N_{x_1}+\theta_2N_{x_2})\right)dx_1dx_2\\
&:=&J_{1,1}+J_{1,2}.
\end{eqnarray*}

Observe that
\begin{eqnarray*}
&&J_{1,1}\\&=&\frac12\textrm{Im}\left[\varepsilon^2\int_{\frac{2s}{\varepsilon^2}}^{\frac{2t}{\varepsilon^2}}\int_{\frac{2s}{\varepsilon^2}}^{x_2}
E\left(\varphi\left(x_{\varepsilon}^{\theta}(s_1),\dots,x_{\varepsilon}^{\theta}(s_k)\right)
e^{i(\theta_1N_{x_1}-\theta_2N_{x_2})}\right)dx_1dx_2\right],
\end{eqnarray*}
that is, $J_{1,1}$ is the imaginary part of the same expression that
$I_{1,1}$ is the real part. So, the same computations for $I_{1,1}$
show the convergence to zero of $J_{1,1}$ when $\varepsilon$ tends
to 0. The same parallelism can be done between the terms $J_{1,2}$,
$J_{2,1}$, $J_{2,2}$ and $I_{1,2}$, $I_{2,1}$, $I_{2,2}$,
respectively.

The proof of Theorem \ref{teo1} is now complete.

\section*{Acknowledgements}
Xavier Bardina and Carles Rovira are partially supported by MEC-FEDER grants MTM2006-06427 and MTM2006-01351, respectively.

\end{document}